\documentclass[12pt,oneside,reqno]{amsart}
\usepackage{txfonts}
\usepackage{bbm}
\usepackage{amsmath}
\usepackage{graphicx}
\usepackage{mathrsfs}
\usepackage{stmaryrd}
\usepackage{amsfonts}
\usepackage{enumerate,amsmath,amssymb,amsthm}
\newcommand{\be}{\begin{eqnarray}}
\newcommand{\ee}{\end{eqnarray}}
\newcommand{\ce}{\begin{eqnarray*}}
\newcommand{\de}{\end{eqnarray*}}
\newtheorem{theorem}{Theorem}[section]
\newtheorem{lemma}[theorem]{Lemma}
\newtheorem{remark}[theorem]{Remark}
\newtheorem{definition}[theorem]{Definition}
\newtheorem{proposition}[theorem]{Proposition}
\newtheorem{Examples}[theorem]{Example}
\newtheorem{corollary}[theorem]{Corollary}

\def\eps{\varepsilon}

\def\p{\partial}

\def\[{{\Big[}}
\def\]{{\Big]}}
\def\<{{\langle}}
\def\>{{\rangle}}
\def\({{\Big(}}
\def\){{\Big)}}

\def\bx{{\mathbf{x}}}

\def\dif{{\mathord{{\rm d}}}}
\def\dis{{\mathord{{\rm \bf d}}}}

\def\min{{\mathord{{\rm min}}}}

\def\no{\nonumber}
\def\={&\!\!=\!\!&}
\def\bt{\begin{theorem}}
\def\et{\end{theorem}}
\def\bl{\begin{lemma}}
\def\el{\end{lemma}}
\def\br{\begin{remark}}
\def\er{\end{remark}}
\def\bd{\begin{definition}}
\def\ed{\end{definition}}
\def\bp{\begin{proposition}}
\def\ep{\end{proposition}}
\def\bc{\begin{corollary}}
\def\ec{\end{corollary}}
\def\bx{\begin{Examples}}
\def\ex{\end{Examples}}

\def\cC{{\mathcal C}}

\def\cF{{\mathcal F}}

\def\cK{{\mathcal K}}

\def\cP{{\mathcal P}}

\def\cW{{\mathcal W}}

\def\mE{{\mathbb E}}

\def\mM{{\mathbb M}}
\def\mN{{\mathbb N}}

\def\mX{{\mathbb X}}
\def\mY{{\mathbb Y}}

\def\sB{{\mathscr B}}
\def\sC{{\mathscr C}}

\def\sF{{\mathscr F}}

\def\sK{{\mathscr K}}

\def\sW{{\mathscr W}}

\def\geq{\geqslant}
\def\leq{\leqslant}

\allowdisplaybreaks

\begin{document}

\title{Stochastic Monge-Kantorovich Problem and its Duality$^*$}

\date{}
\author{Xicheng Zhang}
\thanks{$*$ This work is supported by NSFs of China (Nos. 10971076; 10871215).}
\address{Xicheng Zhang:
Department of Mathematics,
Huazhong University of Science and Technology\\
Wuhan, Hubei 430074, P.R.China\\
email: XichengZhang@gmail.com
}

\begin{abstract}
In this article we prove the existence of a stochastic optimal transference
plan for a stochastic Monge-Kantorovich problem
by measurable selection theorem. A stochastic version of Kantorovich duality and the characterization
of stochastic optimal  transference plan are also established.
Moreover, Wasserstein distance between two probability kernels are  discussed too.
\end{abstract}

\maketitle
\rm

\section{Introduction and Main Results}
Let $\mX$ be a Polish space and $\cP(\mX)$  the total of probability measures on $(\mX,\sB(\mX))$, where
$\sB(\mX)$ is the Borel $\sigma$-field. It is well known that $\cP(\mX)$ is a Polish space with respect to
the weak convergence topology. Let $\sB(\cP(\mX))$ be
the associated Borel $\sigma$-field.  Let $\mY$ be another Polish space and
$c:\mX\times\mY\to[0,\infty]$ be a lower semicontinuous function called cost function.
For $\mu\in\cP(\mX)$ and $\nu\in\cP(\mY)$, consider the classical Monge-Kantorovich problem
\be
C^{\mathrm{deter}}(c,\mu,\nu):=\inf_{\pi\in\Pi(\mu,\nu)}\int_{\mX\times\mY}c(x,y)\pi(\dif x, \dif y),\label{Det}
\ee
where $\Pi(\mu,\nu)$ denotes the set of all joint probability measures on $\mX\times\mY$ with marginal
distributions $\mu$ and $\nu$. The history and the background of Monge-Kantorovich problem are refereed to \cite{Ra-Ru,Vi} etc. The element in $\Pi(\mu,\nu)$ is called transference plan; those achieving
the infimum are called optimal transference plan. We remark that
the existence of optimal transference plan is easily obtained
by the compactness of $\Pi(\mu,\nu)$ in $\cP(\mX\times\mY)$. Moreover, the following Kantorovich duality formula holds
(cf. \cite{Ra-Ru} or \cite[Theorem 5.10]{Vi})
\be
C^{\mathrm{deter}}(c,\mu,\nu)=\sup_{(\psi,\phi)\in L^1(\mu)\times L^1(\nu); \phi-\psi\leq c}
\left(\int_\mY\phi(y)\nu(\dif y)-\int_\mX\psi(x)\mu(\dif x)\right).
\ee

We now turn to the description of stochastic versions of Monge-Kantorovich problem and its duality.
Let $(\Omega,\sF, P)$ be a probability space and $\mu$ a probability kernel from $\Omega$ to $\mX$.
Here, by a probability kernel $\mu$ from $\Omega$ to $\mX$, we mean that a mapping
$\mu:\Omega\times\sB(\mX)\to[0,1]$ satisfies
$$
\mbox{(i) for each $\omega\in\Omega$,
$\mu_\omega\in\cP(\mX)$; (ii) for each $B\in\sB(\mX)$, $\omega\mapsto\mu_\omega(B)$ is $\sF$-measurable.}
$$
Let $\mY$ be another Polish space and $\nu$ a  probability kernel from $\Omega$ to $\mY$.
Let $c:\Omega\times\mX\times\mY\to[0,\infty]$ be a measurable function called
stochastic cost function. Consider
the following stochastic Monge-Kantorovich problem:
\be
C^{\mathrm{stoch}}(c,\mu,\nu):=\inf_{\pi\in\cK(\mu,\nu)}\mE\int_{\mX\times\mY}c(\omega,x,y)\pi_\omega(\dif x,\dif y),
\label{Sto}
\ee
where $\cK(\mu,\nu)$ is the set of all probability kernels from $\Omega$ to $\mX\times\mY$
with marginal probability kernels $\mu$ and $\nu$, i.e., for a $\pi_\omega\in\cK(\mu,\nu)$,
$$
\pi_\omega(\cdot,\mY)=\mu_\omega,\ \
\pi_\omega(\mX,\cdot)=\nu_\omega.
$$
If $\pi^{\mathrm{opt}}\in\cK(\mu,\nu)$ attains the infimum for the minimization problem (\ref{Sto}),
we call it a stochastic optimal transference plan. Unlike the deterministic problem (\ref{Det}), it seems to be hard to
prove the existence of a stochastic optimal transference plan by a direct compactness argument.
In fact, when the cost function $c$ is deterministic, the existence of $\pi^{\mathrm{opt}}_\omega$ has been obtained
by Zhang \cite{Zh} (see also \cite[Corollary 5.22]{Vi}).
On the other hand, one may also expect the following stochastic Kantorovich duality formula holds:
\be
C^{\mathrm{stoch}}(c,\mu,\nu)=\sup_{(\psi,\phi)\in L^1(\mu_\omega\times P)\times L^1(\nu_\omega\times P); \phi-\psi\leq c}
\mE\left(\int_\mY\phi(\omega, y)\nu_\omega(\dif y)-\int_\mX\psi(\omega, x)\mu_\omega(\dif x)\right),
\ee
where $L^1(\mu_\omega\times P)$ denotes the set of all measurable functions $\psi$ with
$\mE\int_\mX|\psi(\omega,x)|\mu_\omega(\dif x)<+\infty$, and $\phi-\psi\leq c$ means that
$\phi(\omega,y)-\psi(\omega,x)\leq c(\omega,x,y)$ for all $\omega,x,y$.

Our first result is about the existence of stochastic optimal transference plans.
\bt\label{Main}
Assume that for each $\omega$, $(x,y)\mapsto c(\omega,x,y)$ is continuous, and for each $(x,y)\in\mX\times\mY$,
$\omega\mapsto c(\omega,x,y)$ is $\sF$-measurable and satisfies
\be
\mE\int_{\mX\times\mY}  c(\omega,x,y)\mu_\omega(\dif x)\nu_\omega(\dif y)<+\infty.\label{Op4}
\ee
Then there exists a stochastic optimal transference plan $\pi^{\mathrm{opt}}\in\cK(\mu,\nu)$ such that
\be
C^{\mathrm{stoch}}(c,\mu,\nu)=\mE\int_{\mX\times\mY}c(\omega,x,y)\pi^{\mathrm{opt}}_\omega(\dif x,\dif y)<+\infty.\label{Op3}
\ee
Moreover, $\omega\mapsto C^{\mathrm{deter}}(c(\omega),\mu_\omega,\nu_\omega)$ is $\sF$-measurable and we have
\be
C^{\mathrm{stoch}}(c,\mu,\nu)
=\mE\left(\inf_{\pi\in\Pi(\mu_\omega,\nu_\omega)}\int_{\mX\times\mY}c(\omega,x,y)\pi(\dif x,\dif y)\right)
=\mE\Big( C^{\mathrm{deter}}(c(\omega),\mu_\omega,\nu_\omega)\Big).\label{PP1}
\ee
\et
\br
For fixed $\omega\in\Omega$, let $X_\omega\subset\Pi(\mu_\omega,\nu_\omega)$ be
the set of all optimal transference plans for deterministic problem (\ref{Det}).
It is well known that $X_\omega$ is a nonempty compact subset of $\cP(\mX\times\mY)$. For proving Theorem \ref{Main},
we have to carefully choose a measurable function $\omega\to\pi^{\mathrm{opt}}_\omega$ so that for each $\omega$,
$\pi^{\mathrm{opt}}_\omega\in X_\omega$. This seems not to be trivial as shown in \cite{Zh}.
\er

Our second result is about the stochastic Kantorovich duality.
\bt\label{Main2}
Keeping the same assumptions as in Theorem \ref{Main}, we further have
\be
C^{\mathrm{stoch}}(c,\mu,\nu)
&=&\sup_{(\psi,\phi)\in L^1(\mu_\omega\times P)\times L^1(\nu_\omega\times P); \phi-\psi\leq c}
\mE\left(\int_\mY\phi(\omega, y)\nu_\omega(\dif y)-\int_\mX\psi(\omega, x)\mu_\omega(\dif x)\right)\no\\
&=&\sup_{(\psi,\phi)\in {Lip}^\omega_b(\mX)\times {Lip}^\omega_b(\mY); \phi-\psi\leq c}
\mE\left(\int_\mY\phi(\omega, y)\nu_\omega(\dif y)-\int_\mX\psi(\omega, x)\mu_\omega(\dif x)\right),
\ee
where ${Lip}^\omega_b(\mX)$ is the space of all bounded measurable functions $\psi(\omega,x)$
on $\Omega\times \mX$ which is Lipschitz continuous in $x$ for each $\omega$, similarly for ${Lip}^\omega_b(\mY)$.
\et

Our third result is about the characterization of stochastic optimal transference plan, which corresponds
to \cite[Theorem 5.10 (ii)]{Vi} (see also \cite{Am, Sc-Te}).
\bt\label{Main3}
In the situation of Theorem \ref{Main}, for any $\pi\in\cK(\mu,\nu)$, the following statements are equivalent:

(a) $\pi$ is a stochastic optimal transference plan;

(b) for almost all $\omega\in\Omega$, the support of $\pi_\omega$ is a $c(\omega)$-cyclically monotone set;

(c) there exist a pair of measurable functions $(\phi,\psi)$ on $\Omega\times\mY$ and
$\Omega\times\mX$ such that
$$
\phi(\omega,y)-\psi(\omega,x)\leq c(\omega,x,y),\ \ \forall (\omega,x,y)\in\Omega\times\mX\times\mY,
$$
and for each $\omega\in\Omega$, $\psi(\omega)$ is $c(\omega)$-convex and
$$
\Gamma_\omega:=\{(x,y):\phi(\omega,y)-\psi(\omega,x)=c(\omega,x,y)\}\subset\p_c\psi(\omega)
$$
has $\pi_\omega$-full measure, where $\p_c\psi(\omega)$ denotes the
$c(\omega)$-subdifferential of $\psi(\omega,\cdot)$.

Moreover, the measurable set $\Gamma:=\{(\omega,x,y): (x,y)\in\Gamma_\omega\}$ defined from (c)
may be  independent  of the choice of optimal plan $\pi$. More precisely, let $\tilde\pi$ be another stochastic
optimal plan,  then $\tilde\pi_\omega$ is concentrated on $\Gamma_\omega$ for almost all $\omega$.
\et
\br
In these theorems, if we assume that $c$ is lower semi-continuous and approximate it by the usual Lipscitz continuous
functions (see (\ref{Ep}) below), then we shall encounter a very subtle issue about the measurability of an uncountable
infimum of lower semi-continuous functions (cf. \cite[p.70-72]{Vi}).
\er

These three theorems will be proved in Section 3 by measurable selection theorem.
For this aim, we give some necessary preliminaries in Section 2. In Section 4, we shall give a definition of Wasserstein distance
between two probability kernels and discuss the corresponding properties. It is hoped that the results of the present paper
can be used to the study of Markov processes.

\section{Preliminaries}

Let $\sC$ be the total of all nonnegative continuous cost functions $c:\mX\times\mY\to[0,\infty)$,
which is endowed with a metric as follows:
$$
\dis_\sC(c_1,c_2):=\sum^\infty_{m=1}2^{-m}\left(1\wedge
\sup_{(x,y)\in B^m_\mX(x_0)\times B^m_\mY(y_0)}|c_1(x,y)-c_2(x,y)|\right),
$$
where $(x_0,y_0)\in\mX\times\mY$ is fixed and
$$
B^m_\mX(x_0):=\{x\in\mX: \dis_\mX(x,x_0)\leq m\},\ \
B^m_\mY(y_0):=\{y\in\mY: \dis_\mY(y,y_0)\leq m\}.
$$
It is easy to see that $(\sC,\dis_\sC)$ is a complete metric space. Let $\mM$ be defined by
$$
\mM:=\left\{(c,\mu,\nu)\in\sC\times\cP(\mX)\times\cP(\mY):
\int_{\mX\times\mY}c(x,y)\mu(\dif x)\nu(\dif y)<+\infty\right\}.
$$
Then it is a metric space (maybe not complete and separable) under
$$
\dis_\mM((c_1,\mu_1,\nu_1),(c_2,\mu_2,\nu_2)):=\dis_\sC(c_1,c_2)+\dis_{\cP(\mX)}(\mu_1,\mu_2)
+\dis_{\cP(\mY)}(\nu_1,\nu_2),
$$
where $\dis_{\cP(\mX)}$ and $\dis_{\cP(\mY)}$ are weak convergence metric in $\cP(\mX)$ and $\cP(\mY)$ respectively.
We have:

\bl\label{Le3}
Let $\{(c_n,\mu_n,\nu_n)\in\mM, n\in\mN\}$ satisfy that
$$
\sup_{n\in\mN}\int_{\mX\times\mY}c_n(x,y)\mu_n(\dif x)\nu_n(\dif y)\leq M.
$$
Assume that $(c_n,\mu_n,\nu_n)$ converges to $(c,\mu,\nu)$ in $\mM$. Then
$$
\int_{\mX\times\mY}c(x,y)\mu(\dif x)\nu(\dif y)\leq M.
$$
\el
\begin{proof}
By Urysohn's lemma, there exist continuous functions $f^m_\mX:\mX\to[0,1]$ and $f^m_\mY:\mY\to[0,1]$ such that
$$
f^m_\mX(x)=1, \ \ x\in B^m_\mX(x_0),\ \  f^m_\mX(x)=0, \ \ x\notin B^{m+1}_\mX(x_0)
$$
and
$$
f^m_\mY(y)=1, \ \ y\in B^m_\mY(y_0),\ \  f^m_\mY(y)=0, \ \ y\notin B^{m+1}_\mY(y_0).
$$
Thus, by the monotone convergence theorem, we have
\ce
\int_{\mX\times\mY}c(x,y)\mu(\dif x)\nu(\dif y)
&=&\lim_{m\to\infty}\int_{\mX\times\mY}c(x,y)\wedge m\cdot f^m_\mX(x)f^m_\mY(y)\mu(\dif x)\nu(\dif y)\\
&=&\lim_{m\to\infty}\lim_{n\to\infty}\int_{\mX\times\mY}c(x,y)\wedge m\cdot f^m_\mX(x)f^m_\mY(y)\mu_n(\dif x)\nu_n(\dif y).
\de
Since $c_n\to c$ in $\sC$, we have
$$
\lim_{n\to\infty}\sup_{(x,y)\in B^{m+1}_\mX(x_0)\times B^{m+1}_\mY(y_0)}|c(x,y)-c_n(x,y)|=0.
$$
Hence,
$$
\int_{\mX\times\mY}c(x,y)\mu(\dif x)\nu(\dif y)
=\lim_{m\to\infty}\lim_{n\to\infty}\int_{\mX\times\mY}c_n(x,y)\wedge m
\cdot f^m_\mX(x)f^m_\mY(y)\mu_n(\dif x)\nu_n(\dif y)\leq M.
$$
The proof is complete.
\end{proof}

We recall the following definitions of cyclical monotonicity and $c$-convexity
(cf. \cite[Definitions 5.1, 5.2]{Vi}).
\bd
Let $\mX,\mY$ be two arbitrary set and $c:\mX\times\mY\to(-\infty,\infty]$ be a function. A subset
$\Gamma\subset\mX\times\mY$ is said to be $c$-cyclically monotone if for any $N\in\mN$ and any family
$(x_1,y_1),\cdots,(x_N,y_N)$ of points in $\Gamma$, the following inequality holds:
$$
\sum_{i=1}^N c(x_i,y_i)\leq\sum_{i=1}^N c(x_i,y_{i+1}),\ \ y_{N+1}=y_1.
$$
A function $\psi:\mX\to(-\infty,+\infty]$ is said to be $c$-convex if it is not identically $+\infty$,
and there exists $\zeta: \mY\to[-\infty,+\infty]$ such that
$$
\psi(x)=\sup_{y\in\mY}(\zeta(y)-c(x,y)),\ \forall x\in\mX.
$$
Then its $c$-transform is defined by
$$
\psi^c(y):=\inf_{x\in\mX}(\psi(x)+c(x,y)), \  \forall y\in\mY,
$$
and its $c$-subdifferential  defined by
$$
\p_c\psi:=\{(x,y)\in\mX\times\mY: \psi^c(y)-\psi(x)=c(x,y)\}
$$
is a $c$-cyclically monotone set.
\ed

We first prove the following slight extension of \cite[Theorem 3]{Sc-Te} and \cite[Theorem 5.20]{Vi}.
\bt\label{Th2}
Assume that $(c_n,\mu_n,\nu_n)\to (c,\mu,\nu)$ in $\mM$. Let $\pi_n$ be an optimal transference plan for problem (\ref{Det})
associated with $c_n,\mu_n,\nu_n$. Then there exists a subsequence still denoted by $n$
such that $\pi_n$ weakly converges to some $\pi\in\Pi(\mu,\nu)$ and $\pi$ is an
optimal transference plan associated with $c,\mu, \nu$.
\et
\begin{proof}
First of all, by \cite[Lemma 4.4]{Vi}, $(\pi_n)_{n\in\mN}$ is tight, and so there exists a subsequence still denoted by $n$ weakly converging to
some $\pi\in\Pi(\mu,\nu)$.

By \cite[Theorem 5.10]{Vi}, $\pi_n$ is concentrated on some $c_n$-cyclically monotone set $\Gamma_n$.
For $N\in\mN$, let $\cC_n(N)\subset (\mX\times\mY)^{\otimes N}$ be defined by
$$
\sum_{i=1}^Nc_n(x_i,y_i)\leq\sum_{i=1}^Nc_n(x_i,y_{i+1}),\ \ y_{N+1}=y_1,
$$
where $(x_i,y_i)^N_{i=1}\in(\mX\times\mY)^{\otimes N}$. Then $\pi^{\otimes N}_n$
is concentrated on $\Gamma_n^{\otimes N}\subset\cC_n(N)$.

For any $\eps\in[0,1]$, let $\cC_\eps(N)\subset (\mX\times\mY)^{\otimes N}$ be defined by
$$
\sum_{i=1}^Nc(x_i,y_i)\leq\sum_{i=1}^Nc(x_i,y_{i+1})+\eps,\ \ y_{N+1}=y_1,
$$
where $(x_i,y_i)^N_{i=1}\in(\mX\times\mY)^{\otimes N}$.
Since $c_n\to c$ in $\sC$, for any $\eps\in(0,1]$ and $N,m\in\mN$, there exists a $n_0\in\mN$ such that for all $n\geq n_0$
$$
\cC_n(N)\cap(B^m_\mX(x_0)\times B^m_\mY(y_0))^{\otimes N}\subset \cC_\eps(N)
\cap(B^m_\mX(x_0)\times B^m_\mY(y_0))^{\otimes N}=:A^m_\eps(N).
$$
Since $c$ is continuous, $A^m_\eps(N)$ is closed. Hence,
$$
\pi^{\otimes N}(A^m_\eps(N))\geq\varlimsup_{n\to\infty}\pi^{\otimes N}_n(A^m_\eps(N))
\geq\varlimsup_{n\to\infty}\pi^{\otimes N}_n(\cC_n(N)\cap(B^m_\mX(x_0)\times B^m_\mY(y_0))^{\otimes N}).
$$
In view that $\pi^{\otimes N}_n$ is concentrated on $\cC_n(N)$, by letting $\eps\downarrow 0$, we further have
\be
\pi^{\otimes N}(A^m_0(N))
\geq\varlimsup_{n\to\infty}[\pi_n(B^m_\mX(x_0)\times B^m_\mY(y_0))]^{N}
\geq\left[1-\varliminf_{n\to\infty}(\mu_n((B^m_\mX(x_0))^c)+\nu_n((B^m_\mY(y_0))^c))\right]^N.\label{Lp1}
\ee
Noticing that $(\mu_n)_{n\in\mN}$ and $(\nu_n)_{n\in\mN}$ are tight, we have
$$
\lim_{m\to\infty}\sup_{n\in\mN}\mu_n((B^m_\mX(x_0))^c)=0,\ \ \lim_{m\to\infty}\sup_{n\in\mN}\nu_n((B^m_\mY(y_0))^c)=0.
$$
Therefore, letting $m\to\infty$ for both sides of (\ref{Lp1}), we obtain that
$$
\pi^{\otimes N}(\cC_0(N))=1,\ \ \forall N\in\mN,
$$
which leads to
$$
(\mbox{support of $\pi$})^{\otimes N}=\mbox{support of $\pi^{\otimes N}$}\subset\cC_0(N),\ \ \forall N\in\mN,
$$
So  the support of $\pi$  is $c$-cyclically monotone.
Since $(c,\mu,\nu)\in\mM$, we have
$$
C^{\mathrm{deter}}(c,\mu,\nu)\leq\int_{\mX\times\mY}c(x,y)\mu(\dif x)\nu(\dif y)<+\infty.
$$
By \cite[Theorem 5.10]{Vi} again, $\pi$ is an optimal transference plan associated with $c,\mu,\nu$.
\end{proof}

The following lemma will be used in the proof of Theorem \ref{Main2}.

\bl\label{Co}
Assume that $C_b(\mX\times\mY)\ni c_n\uparrow c$ in the sense of pointwise. Then
$$
C^{\mathrm{deter}}(c,\mu,\nu)\leq\varliminf_{n\to\infty}C^{\mathrm{deter}}(c_n,\mu,\nu).
$$
\el
\begin{proof}
Without loss of generality, we assume that
$$
\alpha:=\varliminf_{n\to\infty}C^{\mathrm{deter}}(c_n,\mu,\nu)<+\infty.
$$
In particular, there exists a subsequence still denoted by $n$ such that
$$
\lim_{n\to\infty}C^{\mathrm{deter}}(c_n,\mu,\nu)=\alpha.
$$
Let $\pi_n\in\Pi(\mu,\nu)$ be the optimal transference plan associated with $c_n,\mu,\nu$. Since $\Pi(\mu,\nu)$ is weakly compact,
there exists another subsequence $n_k$ such that $\pi_{n_k}$ weakly converges to some $\pi_0\in\Pi(\mu,\nu)$.
By the monotonicity of $c_n$, we have for each $m\in\mN$,
\ce
\int_{\mX\times\mY}c_m(x,y)\pi_0(\dif x,\dif y)&=&\lim_{k\to\infty}\int_{\mX\times\mY}c_m(x,y)\pi_{n_k}(\dif x,\dif y)\\
&\leq&\varlimsup_{k\to\infty}\int_{\mX\times\mY}c_{n_k}(x,y)\pi_{n_k}(\dif x,\dif y)\\
&=&\varlimsup_{k\to\infty}C^{\mathrm{deter}}(c_{n_k},\mu,\nu)=\alpha.
\de
On the other hand, by the monotone convergence theorem, we have
$$
C^{\mathrm{deter}}(c,\mu,\nu)\leq\int_{\mX\times\mY}c(x,y)\pi_0(\dif x,\dif y)=\lim_{m\to\infty}
\int_{\mX\times\mY}c_m(x,y)\pi_0(\dif x,\dif y).
$$
The result now follows.
\end{proof}

We also recall the following measurability theorem for multifunctions (cf. \cite{Ca-Va} or \cite[p.26, Theorem 2.3]{Mo}).
\bt\label{Th1}
Let $(W,\sW)$ be a measurable space and $\mX$ a Polish space.
Let $X: W\to \cF$ be a multifunctions, where $\cF$ is the total of all closed sets in $\mX$. Consider the
following statements:

(1) for any closed $A\subset\mX$.
$$
\{w: X(w)\cap A\not=\emptyset\}\in\sW;
$$

(2) for any open set $A\subset\mX$
$$
\{w: X(w)\cap A\not=\emptyset\}\in\sW;
$$

(3) there exists a sequence $(\xi_n)_{n\in\mN}$ of measurable selections of $X$ such that for each $w\in W$
$$
X(w)=\overline{\{\xi_n(w), n\in\mN\}}.
$$

Then it holds that (1)$\Rightarrow$(2)$\Leftrightarrow$(3).
\et

The following lemma is useful.
\bl\label{Le2}
The Borel $\sigma$-field $\sB(\cP(\mX))$ coincides with the $\sigma$-field generated by the mapping
$\mu\mapsto\mu(B)$, where $B\in\sB(\mX)$.
\el
\begin{proof}
Let $F$ be a closed set in $\mX$. Define
$$
f_n(x):=\frac{1}{(1+\dis_\mX(x,F))^n}.
$$
Then $f_n(x)\downarrow 1_F(x)$. So, for any $r\in[0,1]$
$$
\{\mu\in\cP(\mX): \mu(F)<r\}=\cup_{n\in\mN}\{\mu\in\cP(\mX): \mu(f_n)<r\}\in\sB(\cP(\mX)).
$$
The result now follows by a monotone class argument.
\end{proof}

\section{Proofs of Main Theorems}

In this section we give the proofs of Theorems \ref{Main}, \ref{Main2} and \ref{Main3}.
First, we prove Theorem \ref{Main}.

\begin{proof}[Proof of Theorem \ref{Main}]
Define a multi-valued map:
$$
\mM\ni(c,\mu,\nu)\mapsto\Phi(c,\mu,\nu)\subset\cP(\mX\times\mY),
$$
where $\Phi(c,\mu,\nu)$ is the total of all optimal transference plan associated with $c,\mu,\nu$.

By Theorem \ref{Th2}, for each $(c,\mu,\nu)\in\mM$, $\Phi(c,\mu,\nu)$ is a nonempty
compact subset of $\cP(\mX\times\mY)$,
and for any closed set $A\subset\cP(\mX\times\mY)$
$$
\{(c,\mu,\nu)\in\mM_m: \Phi(c,\mu,\nu)\cap A\not=\emptyset\} \mbox{ is a closed subset of $\mM$,}
$$
where $\mM_m:=\left\{(c,\mu,\nu)\in\mM: \int_{\mX\times\mY}c(x,y)\mu(\dif x)\nu(\dif y)\leq m\right\}$.
Indeed, let $(c_n,\mu_n,\nu_n)\in\mM_m$ converge to $(c,\mu,\nu)$. By Lemma \ref{Le3}, we have
$(c,\mu,\nu)\in\mM_m$. Let $\pi_n\in\Phi(c_n,\mu_n,\nu_n)$ weakly converge to some $\pi\in\Pi(\mu,\nu)$.
By Theorem \ref{Th2}, $\pi\in\Phi(c,\mu,\nu)$. Since $A$ is closed, $\pi$ also belongs to $A$.

Note that
$$
\{(c,\mu,\nu)\in\mM: \Phi(c,\mu,\nu)\cap A\not=\emptyset\}=\cup_{m\in\mN}
\{(c,\mu,\nu)\in\mM_m: \Phi(c,\mu,\nu)\cap A\not=\emptyset\}.
$$
By Theorem \ref{Th1}, there exists a $\sB(\mM)/\sB(\cP(\mX\times\mY))$-measurable selection
$(c,\mu,\nu)\mapsto\pi(c,\mu,\nu)$ such that for each $(c,\mu,\nu)\in\mM$
$$
\pi(c,\mu,\nu)\in\Phi(c,\mu,\nu)\subset\Pi(\mu,\nu).
$$
We now define
$$
\pi^{\mathrm{opt}}_\omega:=\pi(c(\omega),\mu_\omega,\nu_\omega).
$$
Since $\omega\mapsto (c(\omega),\mu_\omega,\nu_\omega)$ is $\sF/\sB(\mM)$-measurable by Lemma \ref{Le2}, we thus have
\be
\omega\mapsto\pi^{\mathrm{opt}}_\omega\mbox{ is $\sF/\sB(\cP(\mX\times\mY))$-measurable}.\label{Op6}
\ee
In particular,
$$
\omega\mapsto \int_{\mX\times\mY}c(\omega,x,y)
\pi^{\mathrm{opt}}_\omega(\dif x,\dif y)=C^{\mathrm{deter}}(c(\omega),\mu_\omega,\nu_\omega)
$$
is $\sF$-measurable  and
$$
C^{\mathrm{stoch}}(c,\mu,\nu)\leq \mE\Big(C^{\mathrm{deter}}(c(\omega),\mu_\omega,\nu_\omega)\Big).
$$
The opposite inequality is clear. Thus, we complete the proof of (\ref{Op3}) and (\ref{PP1}).
\end{proof}

\vspace{3mm}

We now prove Theorem \ref{Main2}.

\begin{proof}[Proof of Theorem \ref{Main2}]

 We divide the proof into three steps.

{\bf (Step 1):} First of all, for any $\pi\in\cK(\mu,\nu)$, we have
\be
&&\sup_{(\psi,\phi)\in L^1(\mu_\omega\times P)\times L^1(\nu_\omega\times P); \phi-\psi\leq c}
\mE\left(\int_\mY\phi(\omega, y)\nu_\omega(\dif y)-\int_\mX\psi(\omega, x)\mu_\omega(\dif x)\right)\no\\
&&\qquad=\sup_{(\psi,\phi)\in L^1(\mu_\omega\times P)\times L^1(\nu_\omega\times P); \phi-\psi\leq c}
\mE\left(\int_{\mX\times\mY}(\phi(\omega, y)-\psi(\omega, x))\pi_\omega(\dif x,\dif y)\right)\no\\
&&\qquad\qquad\leq\mE\left(\int_{\mX\times\mY}c(\omega,x, y)\pi_\omega(\dif x,\dif y)\right).\label{PP3}
\ee
Thus, we obtain one side inequality:
$$
\sup_{(\psi,\phi)\in L^1(\mu_\omega\times P)\times L^1(\nu_\omega\times P); \phi-\psi\leq c}
\mE\left(\int_\mY\phi(\omega, y)\nu_\omega(\dif y)-\int_\mX\psi(\omega, x)\mu_\omega(\dif x)\right)
\leq C^{\mathrm{stoch}}(c,\mu,\nu).
$$

{\bf (Step 2):} In this step, we assume that $c(\omega,x,y)$ is bounded and Lipschitz continuous in $(x,y)$ for each $\omega$.

Let $\pi^{\mathrm{opt}}_\omega$ be the stochastic optimal transference plan constructed in Theorem \ref{Main}.
Let $\Gamma_\omega$ be the support of $\pi^{\mathrm{opt}}_\omega$, a $c(\omega)$-cyclically monotone set.
Note that for any open set $A\subset\mX\times\mY$,
$$
\{\omega: \Gamma_\omega\cap A\not=\emptyset\}=\{\omega: \pi_\omega(A)>0\}\in\sF.
$$
By Theorem \ref{Th1}, there exists a sequence $(\xi_n(\omega),\eta_n(\omega))_{n\in\mN}$
of measurable selections of $\Gamma_\omega$ such that for each $\omega\in \Omega$
\be
\Gamma_\omega=\overline{\{(\xi_n(\omega),\eta_n(\omega)), n\in\mN\}}.\label{Op7}
\ee
Define for each $(\omega,x)\in\Omega\times\mX$,
\be
\psi(\omega,x)&:=&\sup_{m\in\mN}\sup_{(x_1,y_1),\cdots, (x_m,y_m)\in\Gamma_\omega}\Big\{[c(\omega, \xi_1(\omega),\eta_1(\omega))-c(\omega, x_1,\eta_1(\omega))]\no\\
&&+[c(\omega,x_1,y_1)-c(\omega,x_2,y_1)]+\cdots+[c(\omega,x_m,y_m)-c(\omega,x,y_m)]\Big\}.\label{D1}
\ee
Arguing as in \cite[p.65, Step 3]{Vi}, we know that
$$\psi(\omega,\xi_1(\omega),\eta_1(\omega))=0
$$
and
$$
\mbox{$\psi(\omega)$ is $c(\omega)$-convex.}
$$
Since $c(\omega,x,y)$ is continuous with respect to $(x,y)$, by (\ref{Op7}) we may write
\be
\psi(\omega,x)&=&\sup_{m\in\mN}\sup_{(x_1,y_1),\cdots, (x_m,y_m)\in\{(\xi_n(\omega),\eta_n(\omega)), n\in\mN\}}
\Big\{[c(\omega, \xi_1(\omega),\eta_1(\omega))-c(\omega, x_1,\eta_1(\omega))]\no\\
&&+[c(\omega,x_1,y_1)-c(\omega,x_2,y_1)]+\cdots+[c(\omega,x_m,y_m)-c(\omega,x,y_m)]\Big\}.\label{D2}
\ee
Hence, for each $x\in\mX$, $\omega\mapsto\psi(\omega,x)$ is $\sF$-measurable.
Moreover, since $c$ is Lipschitz continuous in $(x,y)$, it is easy to see that for each $\omega\in\Omega$,
$x\mapsto\psi(\omega,x)$ is also Lipschitz continuous.
Let $\psi^c(\omega,y)$ be the $c$-transform of $\psi$ defined by
$$
\psi^c(\omega,y):=\inf_{x\in\mX}\Big(\psi(\omega,x)+c(\omega,x,y)\Big).
$$
Then for each $y\in\mY$, $\omega\mapsto\psi^c(\omega,y)$ is also $\sF$-measurable, and
for each $\omega\in\Omega$, $y\mapsto\psi^c(\omega,y)$ is Lipschitz continuous.
Since $c$ is bounded, as in \cite[p.66, Step 4]{Vi}, $\psi^c$ and $\psi$ are bounded.
Note that (cf. \cite[p.65, Step 3]{Vi})
\be
\psi^c(\omega,y)-\psi(\omega,x)=c(\omega,x,y)\mbox{  on $\Gamma_\omega$}.\label{PP5}
\ee
So
$$
\int_\mX\psi^c(\omega,y)\nu_\omega(\dif y)-\int_\mY\psi(\omega,x)\mu_\omega(\dif x)
=\int_{\mX\times\mY}c(\omega,x,y)\pi^{\mathrm{opt}}_\omega(\dif x,\dif y),
$$
which then gives that
$$
C^{\mathrm{stoch}}(c,\mu,\nu)=\mE\left(\int_\mX\psi^c(\omega,y)\nu_\omega(\dif y)-\int_\mY\psi(\omega,x)\mu_\omega(\dif x)\right).
$$

{\bf (Step 3):} For general $c(\omega,x,y)$, define for $n\in\mN$
\begin{align}
c_n(\omega,x,y):=\inf_{(x',y')\in\mX\times\mY}\Big\{\min(c(\omega,x',y'),n)+n\big[\dis_\mX(x,x')
+\dis_\mY(y,y')\big]\Big\}.\label{Ep}
\end{align}
It is easy to see that $c_n$ is Lipschitz continuous, and
$$
c_n(\omega,x,y)\leq\min(c(\omega,x,y),n)
$$
and for each $(\omega,x,y)\in\Omega\times\mX\times\mY$
$$
c_n(\omega,x,y)\uparrow c(\omega,x,y)\ \ n\to\infty.
$$
Thus, by (\ref{PP1}), Lemma \ref{Co} and Fatou's lemma, we have
\be
C^{\mathrm{stoch}}(c,\mu,\nu)&=&\mE\Big( C^{\mathrm{deter}}(c(\omega),\mu_\omega,\nu_\omega)\Big)\no\\
&\leq&\mE\left(\varliminf_{n\to\infty} C^{\mathrm{deter}}(c_n(\omega),\mu_\omega,\nu_\omega)\right)\no\\
&\leq&\varliminf_{n\to\infty}\mE\left( C^{\mathrm{deter}}(c_n(\omega),\mu_\omega,\nu_\omega)\right)\no\\
&=&\varliminf_{n\to\infty}\mE\left(\int_\mX\phi_n(\omega,y)\nu_\omega(\dif y)
-\int_\mY\psi_n(\omega,x)\mu_\omega(\dif x)\right),\label{PP2}
\ee
where $\phi_n=\psi^c_n\in Lip^\omega_b(\mY)$ and $\psi_n\in Lip^\omega_b(\mX)$ constructed in Step 2
satisfy
\be
\phi_n(\omega,y)-\psi_n(\omega,x)\leq c_n(\omega,x,y)\leq c(\omega,x,y).\label{PP4}
\ee
The proof is thus complete by combining with Step 1.
\end{proof}

Lastly, we prove Theorem \ref{Main3}.

\begin{proof}[Proof of Theorem \ref{Main3}]

(a)$\Rightarrow$(b): Let $\pi\in\cK(\mu,\nu)$ be a stochastic optimal transference plan,
and let $(\phi_n,\psi_n)_{n\in\mN}$ be as in (\ref{PP2}). By (\ref{PP3}) and (\ref{PP2}),
we have
$$
\lim_{n\to\infty}
\mE\left(\int_{\mX\times\mY}[c(\omega,x,y)-\phi_n(\omega,y)+\psi_n(\omega,x)]\pi_\omega(\dif x,\dif y)\right)=0.
$$
If necessary, by extracting a subsequence and by (\ref{PP4}),
there is an $\Omega_0\in\sF$ with $P(\Omega_0)=1$ such that
for each $\omega\in\Omega_0$,
$$
\lim_{n\to\infty}
\int_{\mX\times\mY}[c(\omega,x,y)-\phi_n(\omega,y)+\psi_n(\omega,x)]\pi_\omega(\dif x,\dif y)=0.
$$
Fix such an $\omega$. Up to choosing a subsequence (possibly depending on $\omega$), we can assume that for
$\pi_\omega$-almost all $(x,y)\in\mX\times\mY$,
$$
\lim_{n\to\infty}\phi_n(\omega,y)-\psi_n(\omega,x)=c(\omega,x,y).
$$
For $N\in\mN$, by passing to the limit in the inequality
$$
\sum_{i=1}^Nc(\omega,x_i,y_{i+1})\geq\sum_{i=1}^N[\phi_n(\omega,y_{i+1})-\psi_n(\omega,x_i)]
=\sum_{i=1}^N[\phi_n(\omega,y_{i})-\psi_n(\omega,x_i)],
$$
we find that $\pi^{\otimes N}_\omega$ is concentrated on the closed set
$$
\cC_\omega(N):=\left\{(x_i,y_i)_{i=1}^N\in(\mX\times\mY)^{\otimes N}:
\sum_{i=1}^Nc(\omega,x_i,y_{i+1})\geq\sum_{i=1}^Nc(\omega,x_i,y_i)\right\}.
$$
So the support of $\pi_\omega$ is $c(\omega)$-cyclically monotone.

(b)$\Rightarrow$(c): Fix $\pi\in\cK(\mu,\nu)$ and set $\hat\Gamma_\omega:=\mathrm{supp}(\pi_\omega)$.
Since we can redefine $\pi$ on a $P$-negligible set, without loss of generality, we can
assume that for all $\omega\in\Omega$, $\hat\Gamma_\omega$ is $c(\omega)$-cyclically monotone. Define a
$c(\omega)$-convex function $\psi(\omega,x)$ as in (\ref{D1}) in terms of $\hat\Gamma_\omega$.
From (\ref{D2}), we know that $\psi$ is an $\sF\times\sB(\mX)$-measurable function and
for each $\omega$, $x\mapsto\psi(\omega,x)$ is lower semicontinuous. Let $\psi^c(\omega)$
be the $c(\omega)$-transform of $\psi(\omega)$, i.e.,
$$
\psi^c(\omega,y):=\inf_{x\in\mX}\Big(\psi(\omega,x)+c(\omega,x,y)\Big).
$$
Since $\psi^c$ is the infimum of uncountably many measurable functions, it is not known whether $\psi^c$
is $\sF\times\sB(\mY)$-measurable. As in \cite[p.133, Step 2]{Am} or \cite[p.72]{Vi}, we can modify $\psi^c$
on a $\nu_\omega(\dif y)P(\dif \omega)$-negligible set so that it becomes measurable. First, we disintegrate
$\pi_\omega(\dif x,\dif y)P(\dif\omega)$ as $\pi_\omega(\dif x|y)\nu_\omega(\dif y)P(\dif \omega)$ and
define an $\sF\times\sB(\mY)$-measurable function
$$
\hat\phi(\omega,y):=\int_\mX[\psi(\omega,x)+c(\omega,x,y)]\cdot 1_{\hat\Gamma_\omega}(x,y)\pi_\omega(\dif x|y).
$$
Since $\pi_\omega(\hat\Gamma_\omega)=1$ and
$\hat\Gamma_\omega\subset\p_c\psi(\omega)$ (see (\ref{PP5})), there exists
a measurable set $A\in\sF\times\sB(\mY)$ with $\int_A\nu_\omega(\dif y)P(\dif\omega)=1$ such that
for all $(\omega,y)\in A$,
$$
\hat\phi(\omega,y)=\psi^c(\omega,y)\int_\mX1_{\hat\Gamma_\omega}(x,y)\pi_\omega(\dif x|y)=\psi^c(\omega,y).
$$
Let us define an $\sF\times\sB(\mY)$-measurable function by
$$
\phi(\omega,y):=
\left\{
\begin{aligned}
&\hat\phi(\omega,y)=\psi^c(\omega,y),& (\omega,y)\in A;\\
&-\infty,&(\omega,y)\notin A.
\end{aligned}
\right.
$$
Then, it is easy to check that $(\phi,\psi)$ has  the desired properties.

(c)$\Rightarrow$(a): Arguing as in \cite[Theorem 2]{Sc-Te} or \cite[p.72, (d)$\Rightarrow$(a)]{Vi},
we can prove it by a truncation argument.

Moreover, let $\tilde\pi$ be another stochastic optimal
plan, as in \cite[p.73, (a)$\Rightarrow$(e)]{Vi}, we can prove that
$$
\mE\int_{\mX\times\mY}[c(\omega,x,y)-\phi(\omega,y)+\psi(\omega,x)]\tilde\pi_\omega(\dif x,\dif y)=0.
$$
Hence, for almost all $\omega$,
$\tilde\pi_\omega$ is concentrated on
$$
\Gamma_\omega:=\{(x,y)\in\mX\times\mY: \phi(\omega,y)-\psi(\omega,x)=c(\omega,x,y)\}.
$$
The whole proof is finished.
\end{proof}

\section{Wasserstein Metric between Two Probability Kernels}

In this section, we define the Wasserstein metric in the space of all probability kernels and discuss its properties.
Let $(\mX,\dis_\mX)$ be a metric space. For $p\geq 1$, let $\sK_p(\mX)$ be the space of all probability kernels
from $\Omega$ to $\mX$ with
$$
\mE\int_\mX\dis_\mX(x,x_0)^p\mu_\omega(\dif x)<+\infty
$$
for some $x_0\in\mX$ (hence for all $x_0\in\mX$).
Let us define for $\mu,\nu\in\sK_p(\mX)$
$$
\cW_p(\mu,\nu):=\left(\inf_{\pi\in\cK(\mu,\nu)}\mE\int_{\mX\times\mX}\dis_\mX(x,y)^p\pi_\omega(\dif x,\dif y)\right)^{1/p},
$$
which is called $p$-Wasserstein distance.
By Theorem \ref{Main}, we have
\be
\cW_p(\mu,\nu)=\left(\mE W_p(\mu_\omega,\nu_\omega)^p\right)^{1/p},\label{Op0}
\ee
where $W_p(\mu_\omega,\nu_\omega)=C^{\mathrm{deter}}(\dis^p_\mX,\mu_\omega,\nu_\omega)^{1/p}$
is the usual Wasserstein distance between probability measures $\mu_\omega$ and $\nu_\omega$.

The following result is a direct consequence of (\ref{Op0}) and \cite[Theorem 6.18]{Vi}.
\bt
Let $(\mX,\dis_\mX)$ be a complete and separable metric space, and $(\Omega,\sF, P)$ a separable probability space.
Then for any $p\geq 1$, $(\sK_p(\mX),\cW_p)$ is also a complete and separable metric space.
\et

We now consider the case of $p=1$. In this case, Wasserstein distance is usually
called Kantorovich-Rubinstein distance. We have:
\bt
For any $\mu,\nu\in\sK_1(\mX)$,
$$
\cW_1(\mu,\nu)=\sup_{\|\psi(\omega)\|_{Lip}\leq 1}\mE\left(\int_\mX\psi(\omega,x)\nu_\omega(\dif x)
-\int_\mX\psi(\omega,x)\mu_\omega(\dif x)\right),
$$
where
$$
\|\psi(\omega)\|_{Lip}:=\sup_{x,x'\in\mX}\frac{|\psi(\omega,x)-\psi(\omega,x')|}{\dis_\mX(x,x')}.
$$
\et
\begin{proof}
By Theorem \ref{Main2}, it only needs to prove that
\be
&&\sup_{(\psi,\phi)\in {Lip}^\omega_b(\mX)\times {Lip}^\omega_b(\mX); \phi-\psi\leq \dis_\mX}
\mE\left(\int_\mY\phi(\omega, y)\nu_\omega(\dif y)-\int_\mX\psi(\omega, x)\mu_\omega(\dif x)\right)\label{Eq1}\\
&&\qquad=\sup_{\|\psi(\omega)\|_{Lip}\leq 1}\mE\left(\int_\mX\psi(\omega,x)\nu_\omega(\dif x)
-\int_\mX\psi(\omega,x)\mu_\omega(\dif x)\right).\label{Eq2}
\ee
Assume that $\phi(\omega,y)-\psi(\omega,x)\leq\dis_\mX(x,y)$. Then
$$
\phi(\omega,y)\leq\inf_{x\in\mX}(\psi(\omega,x)+\dis_\mX(x,y))=:\psi^{\dis}(\omega,y)
$$
and
$$
\psi(\omega,x)\geq\sup_{y\in\mX}(\psi^{\dis}(\omega,y)-\dis_\mX(x,y))=:\psi^{\dis\dis}(\omega,x).
$$
Thus,
$$
(\ref{Eq1})\leq\sup_{\psi\in{Lip}^\omega_b(\mX)}
\mE\left(\int_\mY\psi^{\dis}(\omega, y)\nu_\omega(\dif y)-\int_\mX\psi^{\dis\dis}(\omega, x)\mu_\omega(\dif x)\right).
$$
On the other hand, it is easy to verify
$$
\|\psi^{\dis}(\omega)\|_{Lip}\leq 1,
$$
and so,
$$
\psi^{\dis}(\omega,x)=\psi^{\dis\dis}(\omega,x).
$$
Hence, (\ref{Eq1})$\leq$(\ref{Eq2}). Moreover, (\ref{Eq1})$\geq$(\ref{Eq2}) is obvious. The proof is complete.
\end{proof}

{\bf Acknowledgements:}

The author is very grateful to Professor Fuqing Gao for telling me the existence of reference \cite{Zh}.
This work is supported by NSFs of China (Nos. 10971076; 10871215).

\end{document}